\documentclass[12pt]{article}
\usepackage{amssymb}
\usepackage{amsmath}
\usepackage{latexsym}

\newtheorem{theorem}{Theorem}[section]
\newtheorem{proposition}[theorem]{Proposition}

\begin{document}

\title{Vertex operator algebras and weak Jacobi forms}
\author{Matthew Krauel\thanks{Supported by an NSA summer graduate fellowship} and Geoffrey Mason\thanks{Supported by the NSA and NSF}\\
Department of Mathematics \\
UC Santa Cruz}

\date{}
\maketitle

\abstract
\noindent
Let $V$ be a strongly regular vertex operator algebra. For a state $h \in V_1$ satisfying appropriate integrality conditions, we prove that
the space spanned by the trace functions Tr$_{M}q^{L(0)-c/24}\zeta^{h(0)}$ ($M$ a $V$-module) is a vector-valued weak Jacobi form of weight $0$ and a certain index $\langle h, h \rangle/2$. We discuss refinements and applications of this result when $V$ is holomorphic, in particular we prove that
if $g=e^{h(0)}$ is a finite order automorphism then Tr$_Vq^{L(0)-c/24}g$ is a modular function of weight $0$ on a congruence subgroup of $SL_2(\mathbb{Z})$. 

\section{Introduction}
The theory of $n$-point functions at genus $g=1$ for regular vertex operator algebras (VOA) and their orbifolds was established in \cite{DLM2} and \cite{Z}. In particular,  the \emph{modular-invariance} (in the sense of vector-valued modular forms) of the space of
partition functions of $V$-modules was proved. The purpose of the present paper is to show how the portion of this theory concerned with partition functions can be extended to a setting in which \emph{elliptic modular forms} are replaced by \emph{weak Jacobi forms}. In particular, we prove that appropriately defined trace functions associated with the irreducible modules of a strongly regular
VOA define a vector-valued weak Jacobi form of weight zero and a certain index. See Theorem 1 below for a precise statement. Other parts of the theory ($n$-point functions for $n \geq 1$ and associated differential equations, for example) can also be extended to the context of weak Jacobi forms (cf.\ \cite{GK}), and one can work more generally with \emph{mutivariable} weak Jacobi forms. These generalizations are more complicated than the case of $0$-point functions that we treat here, and will be considered elsewhere.

\medskip
Special cases of our Theorem 1 can be found in the mathematical literature, e.g.\ in the work of Kac-Peterson on affine Lie algebras
\cite{KP}, \cite{K} and in the work of Dong-Liu-Ma \cite{DLMa}, which uses Jacobi forms associated to lattice VOAs to study rigidity questions in the context of elliptic genera. That there is a close relation between elliptic genera and Jacobi forms has been understood for some time (cf.\ the articles in \cite{L}). There is an extensive physical literature discussing these connections (e.g.\ \cite{KYY}, \cite{GK}), and
the `$M_{24}$ moonshine' proposed by Eguchi \textit{et al} \cite{EOT} involves similar ideas.
Our results may help foster further connections between VOAs, Jacobi forms and elliptic genera. 

\medskip
To state our main results we need to review some background and introduce notation that will be in effect throughout the paper.
 We refer to \cite{LL} for basic facts about VOAs.
 We consider \emph{simple} VOAs $V$ that are \emph{strongly regular} in the sense of \cite{DM2}. This
  means that $V$ satisfies the following conditions: 
 \begin{eqnarray*}
&&(a)\ \mbox{$V$ has $L(0)$-grading} \ V = \mathbb{C}\mathbf{1} \oplus V_1 \oplus \hdots, \\
&&(b)\ \mbox{$V_1$ consists of \emph{primary} states, i.e.}  \ L(1)V_1=0,\\
&&(c)\ \mbox{$V$ is \emph{regular}, i.e.\ \emph{rational} and \emph{$C_2$-cofinite}.}
\end{eqnarray*}

A number of basic properties of $V$ flow from these assumptions.
In particular, because $V$ is rational it has only finitely many (isomorphism classes of) irreducible modules
(\cite{DLM3}). We denote them $(M^1,Y^1),  \hdots, (M^r, Y^r)$, identifying $(M^1, Y^1)$ with the VOA $(V, Y)$. We let
$Y^j(a, z) := \sum_n a^j(n)z^{-n-1}$ denote the vertex operator for $a \in V$
 with respect to $M^j\ (1 \leq j \leq r)$. In case $j=1$ we omit the superscript $j$ from the notation.

\medskip
Thanks to properties (a) and (b), a Theorem of Li  \cite{L1} says that there is a \emph{unique} invariant bilinear form $\langle \ , \ \rangle: V \times V \rightarrow \mathbb{C}$ normalized so that
$\langle \mathbf{1}, \mathbf{1} \rangle = -1$. Furthermore, $\langle \ , \ \rangle$ is \emph{symmetric} and 
(because $V$ is simple) \emph{nondegenerate},  and
$a(1)b = \langle a, b \rangle \mathbf{1}$ for $a, b \in V_1$.

\medskip
For $h \in V_1, \tau \in \mathbb{H}$ (complex upper half-plane) and $z \in \mathbb{C}$,
 we set
\begin{eqnarray}\label{Jdef}
J_{j, h}(\tau, z) := \mbox{Tr}_{M^j} q^{L(0)-c/24}\zeta^{h^j(0)} \ \ (1 \leq j \leq r),
\end{eqnarray}
where $q = e^{2\pi i \tau}, \zeta = e^{2\pi i z}$. We again omit $j$ from the notation when
$j=1$. (\ref{Jdef}) should generally be considered 
as a \emph{formal} sum. However, we will see (cf.\ \cite{DLMa}) that it converges for appropriate choices of $h$.

\medskip
We can now state our main results. 

\medskip
\noindent
{\bf Theorem 1} Suppose that $V$ is a strongly regular simple VOA. Let
$h \in V_1$, and assume that $h^j(0)$ is semisimple
with eigenvalues in $\mathbb{Z}$ for $1 \leq j \leq r$. 
Then $J_{j, h}(\tau, z)$ is holomorphic in $\mathbb{H}\times\mathbb{C}$, and the following functional equations hold for 
all $\gamma=\left(\begin{array}{cc}a&b \\c&d\end{array}\right) \in 
\Gamma,
(u, v)\in \mathbb{Z}^{2}$, and $1 \leq i \leq r$:
\begin{eqnarray*}
&&\ (i)\ \mbox{there are scalars $a_{ij}(\gamma)$ depending only on $\gamma$ such that}\\
&&\ \ \ \ J_{i, h} \left(\gamma\tau, \frac{z}{c\tau+d}\right) =
e^{\pi i cz^2\langle h, h \rangle/(c\tau+d)}\sum_{j=1}^r a_{ij}(\gamma)
J_{j, h}(\tau, z), \\
&&(ii)\ \mbox{there is a permutation $i \mapsto i'$ of $\{1, \hdots, r\}$ such that}\\
&&\ \ \ \ \ J_{i, h}(\tau, z+u\tau+v) 
=e^{-\pi i \langle h, h \rangle(u^2 \tau +2uz)}
J_{i', h}(\tau, z).
\end{eqnarray*}

The sharpest results obtain when $V$ is assumed to be \emph{holomorphic}. That is,
$r=1$ and $V$ is the unique irreducible $V$-module. 

\medskip
\noindent
{\bf Theorem 2} Suppose that $V$ is a holomorphic, strongly regular VOA. Let $h \in V_1$, and assume that $h(0)$ has integral eigenvalues. Then $m := \langle h, h \rangle/2$ is
an integer, and $J_{h}(\tau, z)$ is holomorphic in $\mathbb{H}\times\mathbb{C}$ and satisfies the following functional equations:
\begin{eqnarray*}
&&\ (i)\ J_{h} \left(\gamma\tau, \frac{z}{c\tau+d}\right) =
e^{2\pi i cz^2 m/(c\tau+d)}\chi(\gamma)J_{h}(\tau, z), \\
&&(ii)\  J_{h}(\tau, z+u\tau+v) 
=e^{-2\pi i m(u^2 \tau +2uz)} J_{h}(\tau, z),
\end{eqnarray*}
(\mbox{$\chi: \Gamma \rightarrow \mathbb{C}^*$ is a character}).

\bigskip
\noindent
{\bf Remarks} \\
1. Theorem 1 essentially
says that the column vector of functions 
\begin{eqnarray}\label{colJ}
J(\tau, z) := (J_{1, h}(\tau, z), \hdots, J_{r, h}(\tau, z))^t
\end{eqnarray}
is a \emph{vector-valued weak Jacobi form} of weight $0$ and index $\langle h, h \rangle/2$. This is not quite true as it stands, because $\langle h, h \rangle/2$ is not necessarily an integer unless 
$V$ is holomorphic. However, we will see (Proposition \ref{proprational}) that $\langle h, h \rangle$ is always \emph{rational}, in which case
a more precise statement is that $J(\tau, z)$ is a vector-valued weak Jacobi form on a congruence subgroup of the Jacobi group $\Gamma \ltimes \mathbb{Z}^2$ (i.e. on a subgroup of finite index  that contains a congruence subgroup of $\Gamma$.)\\
2. Theorem 2 says that 
$J_h(\tau, z)$ is a \emph{weak Jacobi form} on $\Gamma$ of weight $0$, index $m$, and character $\chi$. \\
3. Multiplying by $\eta(\tau)^{c/24}$ ($\eta(\tau)$ is the Dedekind eta-function) 
eliminates the character $\chi$ in Theorem 2, and we obtain a holomorphic weak Jacobi form of weight 
$c/2$ (an integer divisible by $4$ because $V$ is holomorphic (\cite{Z})) and index $m$ on the full Jacobi group. So there is a Fourier-Jacobi expansion
\begin{eqnarray}\label{FJexp}
\eta(\tau)^{c/24}J_h(\tau, z) = \sum_{n, r \in \mathbb{Z}} c(n, r)q^n\zeta^r
\end{eqnarray}
where $c(n, r)=0$ unless $n \geq 0$ and 
$r^2 \leq m^2+4mn$ (\cite{EZ}, P. 104 \textit{et seq}).\\
4. It is known (cf.\ the proof of Theorem 3 below) that the character $\chi$ intervening in Theorem 2 has order dividing $3$, which amounts to the assertion that
$\chi(S)=1$. \\
5. There are many elements $h \in V_1$ that satisfy the assumptions of Theorems 1 or 2 - indeed, they
\emph{span} $V_1$.

\bigskip
A (true) Jacobi form is a weak Jacobi form satisfying the stronger conditions that 
(in the notation of (\ref{FJexp})) $c(n, r)=0$ unless $r^2 \leq 4mn$ (\cite{EZ}).
We then have the following 

\medskip
\noindent
{\bf Supplement to Theorem 2} Let the notation and assumptions be as in Theorem 2. 
If $m \leq 4$ then $\eta(\tau)^{c/2}J_h(\tau, z)$ is a Jacobi form of weight $c/2$ and index $m$
on the full Jacobi group.

\medskip
\noindent
We do not know if the Supplement to Theorem 2 holds for \emph{all} $m$.

\bigskip
Theorems 1 and 2 can be used to establish (elliptic) modular-invariance results for various trace functions associated to $V$. Our final Theorem illustrates this idea. If $V$ is holomorphic it is known
(\cite{DLM2}, \cite{DLM3}) that there is (up to isomorphsm) a \emph{unique} $g$-twisted sector $V(g)$
whenever $g \in$ Aut$(V)$ has finite order. We are interested in the case when such an automorphism $g$ lies in the subgroup $L\subseteq$ Aut$(V)$  generated by exponentials $e^{a(0)}\ (a \in V_1)$. (In this special case, there is another proof of existence and uniqueness of 
$V(g)$ in \cite{L2}.)

\medskip
\noindent
{\bf Theorem 3}   Suppose that $V$ is holomorphic and strongly regular, and assume that
$g \in L$ has \emph{finite order}. Then the following hold.
\begin{eqnarray*}
&&\ (i)\ Z_V(g, \tau) := \mbox{Tr}_V g\ q^{L(0)-c/24}\ \mbox{is a modular function of weight $0$}\\
&&\ \ \ \ \ \ \mbox{on a congruence subgroup of $\Gamma$}.\\
&&(ii)\ Z_V(g, S\tau) = \mbox{Tr}_{V(g)} q^{L(0)-c/24}.
\end{eqnarray*}

\medskip
The modular-invariance of trace functions for finite order automorphisms $g$ of rational VOAs is an important  conjecture that is open even in the special case when $V$ is holomorphic. Theorem 3 proves the conjecture for holomorphic $V$ and $g\in L$.  (Conjecturally, $L$ has finite index in Aut$(V)$. For further discussion, see \cite{Ma}). Part (ii) says that the $S$-transform of the $g$-trace $Z_V(g, \tau)$ is \emph{precisely}
the graded dimension of the $g$-twisted sector.  Both parts of the Theorem are commonly assumed in the physics literature.

\medskip
An important tool in the proofs of Theorems 1-3 is a result of Li (\cite{L2}) giving a canonical construction of 
modules and twisted sectors for $V$ when the relevant (finite order) automorphism arises as an exponential
$e^{a(0)}$ with $a \in V_1$. As we will see, this dovetails perfectly with the 
theory of Jacobi forms.

\medskip
We use the additional notation $\Gamma = SL_2(\mathbb{Z}), \ 
S = \left(\begin{array}{cc}0 & -1 \\ 1 & 0\end{array}\right)$.

\section{Proofs}
We keep the notation introduced in Section 1.
 First we give the proof of Theorem 1. With the assumptions of the Theorem, the holomorphy
 of $J_{j, h}(\tau, z)$ is proved in \cite{DLMa}, Proposition 1.8.  Turning to part (i), we utilize a modular-invariance result\footnote{The correct formulation of Miyamoto's Main Theorem
is display (27) of \cite{M}. Display (9) (loc. cit.) contains a typo.}  of Miyamoto \cite{M}. 
To describe this, for states $u, v \in V_1$ and $1\leq j \leq r$, introduce the function
\begin{equation}\label{eq: Phidef}
\Phi_{j}(u, v, \tau) = \mbox{Tr}_{M^j} e^{2 \pi i(v^j(0) + \langle u, v \rangle/2)}
q^{L(0) + u^j(0) +  \langle u, u  \rangle/2 - c/24}.
\end{equation}
Comparing our notation to that of Miyamoto (loc.\ cit.) we use $\Phi$ for what Miyamoto calls
$Z$ and we exchange the r\^{o}les of $u$ and $v$.
 Miyamoto's result can then be stated as follows: there are scalars $a_{ij}(\gamma)\ (1 \leq i, j \leq r)$ independent of $(u, v, \tau)$ such that
\begin{eqnarray}\label{recur}
\Phi_i (u, v, \gamma \tau) = \sum_{j=1}^r a_{ij}(\gamma) \Phi_j(au+cv, bu+dv, \tau)
\end{eqnarray}
for all $\gamma = \left(\begin{array}{cc} a&b \\ c&d \end{array}\right) \in \Gamma$.
Note that
\begin{eqnarray*}
\Phi_j (0, zh, \tau) = J_{j, h}(\tau, z).
\end{eqnarray*}
Then we find using (\ref{recur}) that there are scalars $a_{ij}(\gamma)$ such that
\begin{eqnarray*}
&&J_{i, h} \left(\gamma\tau, \frac{z}{c\tau+d}\right) =
\Phi_i\left(0, \frac{zh}{c\tau+d}, \gamma\tau\right) \\
&=&  \sum_{j=1}^r a_{ij}(\gamma) \Phi_j\left(\frac{czh}{c\tau+d}, \frac{dzh}{c\tau+d}, \tau\right) \\
&=& \sum_{j=1}^r a_{ij}(\gamma) \mbox{Tr}_{M^j}
e^{2\pi i(dzh^{j}(0)/(c\tau+d)+cd\langle zh, zh \rangle/2(c\tau+d)^2} \\
&&\ \ \ \ \  q^{L(0)+czh^{j}(0)/(c\tau+d)+c^2\langle zh, zh \rangle /2(c\tau+d)^2-c/24} \\
&=& \sum_{j=1}^r a_{ij}(\gamma) \mbox{Tr}_{M^j}
e^{2\pi i((c\tau+d)zh^{j}(0)/(c\tau+d)+c(c\tau+d)z^2\langle h, h \rangle/2(c\tau+d)^2} q^{L(0)-c/24} \\
&=& \sum_{j=1}^r a_{ij}(\gamma) \mbox{Tr}_{M^j}
e^{2\pi i(zh^{j}(0)+cz^2\langle h, h \rangle/2(c\tau+d))}  q^{L(0)-c/24} \\
&=& e^{\pi i cz^2\langle h, h \rangle/(c\tau+d)}\sum_{j=1}^r a_{ij}(\gamma)
J_{j, h}(\tau, z).
\end{eqnarray*}
This completes the proof of part (i).

\bigskip
We turn to part (ii).
Set
\begin{eqnarray*}
h' := - uh,
\end{eqnarray*}
and note that the zero modes  of $h'$ are semisimple with eigenvalues in $\mathbb{Z}$. 
This permits us to use
a result of Li (\cite{L2}, Proposition 5.4). In the present circumstance it says that for each index $j$ there is an isomorphism of \emph{weak} $V$-modules
\begin{eqnarray}\label{Liiso3}
 (M^{j'}, Y^{j'}_{\Delta_{h'}(z)}) \cong (M^j, Y^j)
\end{eqnarray}
for some $j'$. Furthermore, the map $j \mapsto j'$ is a permutation of the indices
$\{1, \hdots, r\}$.  (Li's result involves irreducible weak $V$-modules, but because $V$ is regular, these are ordinary irreducible $V$-modules (\cite{DLM1}).) 
Recall the meaning of the correspondence
$Y^{j}_{\Delta_{h'}(z)}$ (\cite{L2}):
\begin{eqnarray*}
\Delta_{h'}(z) &:=& z^{h'(0)}\exp\left\{- \sum_{k \geq 1} \frac{h'(k)}{k}(-z)^{-k} \right\},\\
Y^{j'}_{\Delta_{h'}(z)}(v, z) &:=& Y^{j'}(\Delta_{h'}(z)v, z) \ \ (v \in V).
\end{eqnarray*}

\medskip
In order to make use of (\ref{Liiso3}) we need some simple computations which 
we leave to the reader:
\begin{eqnarray*}
h'(0)\omega = h'(2)\omega = 0, \ h'(1)\omega = h'.
\end{eqnarray*}
Then 
\begin{eqnarray*}
\Delta_{h'}(z)\omega &=& \left(z^{h'(0)}\exp\left\{- \sum_{k \geq 1} 
\frac{h'(k)}{k}(-z)^{-k} \right\}\right) \omega \\
&=& \left(1+h'(1)z^{-1} -\frac{1}{2}h'(2)z^{-2}+\frac{1}{2}h'(1)^2z^{-2}\right)\omega \\
&=& \omega + h'z^{-1}+ \frac{1}{2}\langle h', h' \rangle z^{-2},
\end{eqnarray*}
so that
\begin{eqnarray}\label{LDelta}
L_{\Delta_{h'}}(0) = L(0) + h'(0) + \frac{1}{2}u^2\langle h, h, \rangle Id.
\end{eqnarray}

\medskip
Similarly, 
\begin{eqnarray*}
\Delta_{h'}(z)h &=& \left(z^{h'(0)}\exp\left\{ -\sum_{k \geq 1} 
\frac{h'(k)}{k}(-z)^{-k} \right\}\right) h \\
&=& \left(1+ h'(1)z^{-1} -\frac{1}{2}h'(2)z^{-2}+\frac{1}{2}h'(1)^2z^{-2}\right) h \\
&=& h+ \langle h', h \rangle z^{-1},
\end{eqnarray*}
whence
\begin{eqnarray*}
Y_{\Delta_{h'}}(h, z)  = Y(h, z) + \langle h', h \rangle z^{-1}.
\end{eqnarray*}

\medskip
Now we are in a position to use the isomorphism (\ref{Liiso3}). Namely, operators that correspond to each other under the isomorphism necessarily have the \emph{same} trace function.   
The preceding calculations then give
\begin{eqnarray*}
&&J_{j, h}(\tau, z+ u\tau+ v) \\
&=&  \mbox{Tr}_{M^j}q^{L(0)-c/24}e^{2\pi i(z+u\tau)h^j(0)} \\
&=&   \mbox{Tr}_{M^{j'}}q^{L(0)+h'^{j'}(0) +u^2\langle h, h \rangle/2-c/24}
e^{2\pi i (z+u\tau)(h^{j'}(0) +\langle h',  h \rangle)} \\
&=& e^{\pi i u^2\langle h, h \rangle \tau }e^{2\pi i \langle h',  u h\rangle \tau}e^{2\pi iz \langle h', h \rangle}
\mbox{Tr}_{M^{j'}}q^{L(0) - u h^{j'}(0) -c/24} e^{2\pi i (z+u\tau)h^{j'}(0)} \\
&=&e^{-\pi i u^2\langle h, h \rangle \tau }e^{2\pi i z \langle h', h \rangle} \mbox{Tr}_{M^{j'}}q^{L(0)  -c/24} 
e^{2\pi i(z+u\tau)h^{j'}(0)} 
e^{-2\pi i \tau  uh^{j'}(0)} \\
&=&e^{-\pi i \langle h, h \rangle (u^2 \tau +2uz)}
J_{j', h}(\tau, z). 
\end{eqnarray*}
This completes the proof of Theorem 1.

\bigskip
Next we prove the rationality results concerning $\langle h, h \rangle$ discussed in Section 1. Indeed, we prove a more precise result. Before stating this we recall that the $L(0)$-grading on the irreducible
$V$-module $M^j$ has the form
\begin{eqnarray*}
\mbox{Tr}_{M^j}q^{L(0)} = \sum_{n \geq 0} \dim M^j_{n+\lambda_j} q^{n+\lambda_j}
\end{eqnarray*}
for a scalar $\lambda_j$ called the \emph{conformal weight} of $M^j$. Because $V$ is strongly regular,
it is known (\cite{DLM2}) that $\lambda_j \in \mathbb{Q}$. We now have
\begin{proposition}\label{proprational} Assume that $h \in V_1$ satisfies the hypotheses stated in Theorem 1.
Then the following hold:
\begin{eqnarray*}
&&\ (i)\ \mbox{We have $\langle h, h \rangle/2 \equiv \lambda_j(mod\ \mathbb{Z}$) for some index $j$.}\\
&&\ \ \ \ \mbox{ In particular, $\langle h, h \rangle \in \mathbb{Q}$.}\\
&&(ii)\ \mbox{If $V$ is holomorphic then}\ \langle h, h \rangle/2 \in \mathbb{Z}.
\end{eqnarray*}
\end{proposition}
{\bf Proof:} Notice that if $V$ is holomorphic then $\lambda_1=0$ is the only conformal weight. So (ii) is an immediate consequence of (i).

\medskip
For part (i) we use  (\ref{Liiso3}) with $h'$ replaced by $h$. So there is an isomorphism of $V$-modules 
\begin{eqnarray*}
\varphi: (M^{j'}, Y^{j'}_{\Delta_{h}(z)}) \stackrel{\cong}{\rightarrow} (M^j, Y^j)
\end{eqnarray*}
for some $j'$, and
\begin{eqnarray}\label{varconj}
Y^j(v, z)\varphi = \varphi Y^{j'}(\Delta_h(z)v, z) \ \ (v \in V).
\end{eqnarray}
Taking zero modes with $v = \omega$ in (\ref{varconj}) and using calculations essentially identical
to those used in the proof of part (ii) of Theorem 1, we find
\begin{eqnarray*}
\varphi^{-1} L^j(0)\varphi = L^{j'}(0)+h^{j'}(0)+\langle h, h \rangle/2.
\end{eqnarray*}
Choose
$j'=1$, so that $M^{j'}=V$, and apply both sides to the vacuum vector to obtain
\begin{eqnarray*}
\varphi^{-1} L^j(0)\varphi(\mathbf{1}) = \frac{1}{2}\langle h, h \rangle \mathbf{1}.
\end{eqnarray*}
It follows that there is an integer $n_0$ such that $\varphi(\mathbf{1}) \in M^j_{n_0+\lambda^j}$, and
\begin{eqnarray*}
\langle h, h \rangle/2 = n_0 + \lambda^j.
\end{eqnarray*}
This proves (i), and completes the proof of the Proposition. $\hfill \Box$

\bigskip
Using Theorem 1 and Proposition \ref{proprational}, it can be shown that the column vector
(\ref{colJ}) is a vector-valued weak Jacobi form of weight $0$ and index $\langle h, h \rangle/2$ on a congruence subgroup of the Jacobi group. We will not need this result here, so we skip the details.

\medskip
If $V$ is holomorphic, the functional equations (i) and (ii) in Theorem 1 simplify to the corresponding statements in Theorem 2. The constant $a_{11}(\gamma)$, now written $\chi(\gamma)$, is necessarily a character of $\Gamma$. This completes the proof of all of the assertions of Theorem 2, which is now proved. $\hfill \Box$

\bigskip
The supplement to Theorem 2 is a consequence of the following purely arithmetic statement about weak Jacobi forms. 
\begin{proposition}\label{propmleq4} Let $\phi_{k, m}(\tau, z) = \sum_{n, r} c(n, r)q^n\zeta^r$ be a holomorphic weak Jacobi form of weight $k$ and index $m$ on the full Jacobi group. Suppose that $k \geq 4$ and $1 \leq m \leq 4$. Then $\phi_{k, m}$ is a Jacobi form if, and only if, $c(0, r) = 0$ for $r \not=0$.
\end{proposition}
{\bf Proof} We know (\cite{EZ} P.108) that for even $k$ there is a linear isomorphism
\begin{eqnarray*}
&&P: M_k \oplus M_{k+2} \oplus \hdots \oplus M_{k+2m} \rightarrow  \tilde{J}_{k, m}, \\
&& \hspace{2.8cm} (f_0, f_1, \hdots, f_m) \mapsto \sum_{i=0}^m f_i \tilde{\phi}_{-2, 1}^i \tilde{\phi}^{m-i}_{0, 1}.
\end{eqnarray*}
Notation here is as in \cite{EZ}, i.e.\ $M_k$ is the space of holomorphic modular forms
of weight $k$ on $\Gamma$; $\tilde{J}_{k, m}$ the space of holomorphic weak Jacobi forms of weight $k$ and index $m$; and $\tilde{\phi}_{-2, 1},  \tilde{\phi}_{0, 1}$ are weak
Jacobi forms of index $1$ and weights $-2, 0$ respectively, defined by
\begin{eqnarray*}
\tilde{\phi}_{-2, 1} &=& \frac{\phi_{10, 1}}{\Delta} = (\zeta -2 + \zeta^{-1})+ \hdots,  \\
 \tilde{\phi}_{0, 1} &=& \frac{\phi_{12, 1}}{\Delta}  = (\zeta +10 + \zeta^{-1})+ \hdots
\end{eqnarray*} 
Moreover (loc.\ cit.\ Theorem 9.2, Corollary, and Section 10) if $k \geq 3$ then
\begin{eqnarray}\label{codimform}
\dim  J_{k, m} = \dim \tilde{J}_{k, m} - \sum_{\nu=0}^m \left\lceil \frac{\nu^2}{4m} \right \rceil,
\end{eqnarray}
where $J_{k, m} \subseteq \tilde{J}_{k, m}$ is the space of holomorphic Jacobi forms of weight $k$ and index $m$.

\medskip
 Certainly all forms in $J_{k, m}$ satisfy the condition $c(0, r)=0$ for $r\not= 0$.
 We will show that
the space $E$ of functions in $\tilde{J}_{k, m}$ satisfying these vanishing conditions has codimension  exactly $m$.
Consider the weak Jacobi forms $P(E_i), \ 0 \leq i \leq m$, where $E_i = 1+ \hdots$
is an Eisenstein series in $M_{k+2i}$. Then
\begin{eqnarray*}
P(E_i) =  E_i \tilde{\phi}_{-2, 1}^i \tilde{\phi}^{m-i}_{0, 1} = 
(\zeta -2 + \zeta^{-1})^i(\zeta+ 10 + \zeta^{-1})^{m-i} + O(q)
\end{eqnarray*}
 Let $x = \zeta+\zeta^{-1}$. Then
\begin{eqnarray*}
P(E_i) = 
(x -2)^i(x+10)^{m-i} + O(q).
\end{eqnarray*}
One sees inductively that the polynomials $(x -2)^i(x+10)^{m-i}
\ (0 \leq i \leq m)$ are a \emph{basis} for the space of all polynomials in $\mathbb{Q}[x]$ of degree at most $m$.
So there are $m+1$ weak Jacobi forms $Q_i, \ 0 \leq i \leq m$, constructed as certain linear combinations of the $P(E_i)$, such that
\begin{eqnarray*}
Q_i = x^i + O(q),
\end{eqnarray*}
whence $E$  clearly has 
codimension $m$, as asserted.

\medskip
Finally, if $1 \leq m \leq 4$ then 
\begin{eqnarray*}
 \sum_{\nu=0}^m \left\lceil \frac{\nu^2}{4m} \right \rceil = m.
\end{eqnarray*}
Then $J_{k, m}$ has codimension $m$ in $\tilde{J}_{k, m}$ by (\ref{codimform}), and because
$J_{k, m} \subseteq E \subseteq \tilde{J}_{k, m}$ the desired equality $E = J_{k, m}$ holds.
This completes the proof of the Proposition. $\hfill \Box$

\bigskip
Finally, we discuss the proof of Theorem 3. First note that because $V$ is rational then 
it  is finitely generated (\cite{DZ}). It then follows that Aut$(V)$ is an algebraic group (\cite{DG}). The strong regularity of $V$ also implies that $V_1$ is a \emph{reductive} Lie algebra (\cite{DM1}), so that
the subgroup $L \subseteq$ Aut$(V)$ generated by the exponentials $e^{a(0)}\ (a \in V_1)$
is a (normal) reductive algebraic subgroup. Then if $g \in L$ has finite order, say $R$, then
$g=e^{2\pi i a(0)}$ for some semisimple element $a \in V_1$ and $Ra(0)$ has eigenvalues in $\mathbb{Z}$.

\medskip
Now $J_{Ra}(\tau, z)$ is a weak Jacobi form by Theorem 2. By \cite{EZ}, Theorem 1.3 (the proof of which applies to weak Jacobi forms), we find that
$J_{Ra}(\tau, z_0)$ is a modular form on a congruence subgroup of $\Gamma$ whenever
$z_0 \in \mathbb{Q}$. Taking $z_0=R^{-1}$, it follows that
\begin{eqnarray*}
J_{Ra}(\tau, 1/R) = \mbox{Tr}_Vg q^{L(0)-c/24} =: Z_V(g, \tau).
\end{eqnarray*}
is a modular form on a congruence subgroup. This proves part (i) of Theorem 3.

\medskip
Turning to part (ii), we can use part (i) of Theorem 2 to see that
\begin{eqnarray*}
Z_V(g, S\tau) &=& J_{Ra}(S\tau, \tau/R\tau) \\
&=& e^{2\pi i \langle a, a \rangle \tau/2}\chi(S)J_{Ra}(\tau, \tau/R) \\
&=&  \chi(S) \mbox{Tr}_V q^{2\pi i(L(0)+a(0)/R+\langle a, a \rangle/2-c/24)}.
\end{eqnarray*}

\medskip
On the other hand, we can use Li's Theorem (\cite{L2}, Proposition 5.4), this time for
\emph{twisted sectors}. More exactly, there is an isomorphism of twisted sectors
$(V(g), Y) \cong (V, Y_{\Delta_{-a}(z)})$, and by a computation very similar to that leading to
(\ref{LDelta}) we find that
\begin{eqnarray*}
\mbox{Tr}_Vq^{L_{\Delta_{-a}}(0)-c/24}= \mbox{Tr}_V q^{L(0)+a(0)/R+\langle a, a \rangle/2-c/24}.
\end{eqnarray*}
From the preceding two displays it follows that
\begin{eqnarray*}
Z_V(g, S\tau) = \chi(S)\mbox{Tr}_{V(g)}q^{L(0)-c/24}.
\end{eqnarray*}

Therefore, in order to complete the proof of part (ii) of Theorem 3, it suffices to show that
$\chi(S)=1$. To see this, note from Theorem 1 that the scalars $a_{ij}(\gamma)$, and in particular the character
values $\chi(\gamma)$ in the holomorphic case, are independent of $z$. We may therefore set $z=0$, in which case $J_h(\tau, z) = \mbox{Tr}_V q^{L(0)-c/24} =: Z_V(\tau)$ is just the usual partition function for $V$, and
part (i) of Theorem 2  reduces to  the identity $Z_V(S\tau) = \chi(S)Z_V(\tau)$. But it is well-known
that if $V$ is holomorphic then in fact $Z_V(S\tau)=Z_V(\tau)$. (For a proof, see e.g.\ \cite{H}, Korollar 2.1.3.)
So indeed $\chi(S)=1$, and the proof is complete. $\hfill \Box$


\begin{thebibliography}{BPZ}
\bibitem[DG]{DG} C. Dong and R. Griess, Automorphism groups and derivation algebras of finitely generated vertex operator algebras, Mich. Math. J. \textbf{50} (2002), 227-239.


 \bibitem[DLM1]{DLM1} C. Dong, H. Li and G. Mason, Regularity of rational vertex operator algebras,
  Adv. Math. \textbf{132} (1997), 148-166.

\bibitem[DLM2]{DLM2} C. Dong, H. Li and G. Mason, Modular-invariance of trace functions in orbifold theory and generalized moonshine, Comm. Math. Phys. \textbf{214} (2000), 1-56.

\bibitem[DLM3]{DLM3} C. Dong, Li and G. Mason, Twisted representations of vertex operator algebras,
Math. Ann. \textbf{310} (1997), 148-166.

 \bibitem[DLMa]{DLMa} C. Dong, K. Liu, and X. Ma, Elliptic genus and vertex operator algebras, Pure Appl. Math. Quart. \textbf{1} (2005), 791-815.
  

\bibitem[DM1]{DM1} C. Dong and G. Mason, Rational vertex operator algebras and the effective central charge,
Int. Math. Res. Not. \textbf{56} (2004), 2989-3008.

\bibitem[DM2]{DM2} C. Dong and G. Mason, Integrability of $C_2$-cofinite Vertex Operator Algebras,
Int. Math. Res. Not. Vol. 2006, Art. ID \textbf{80468}, 1-15.


\bibitem[DZ]{DZ} C. Dong and W. Zhang, Rational vertex operator algebras are finitely generated, J. Alg. \textbf{320} No. 6
(2008), 2610-2614.

\bibitem[EZ]{EZ} M. Eichler and D. Zagier, \textit{The Theory of Jacobi Forms}, Progress in Math. Vol.
\textbf{55}, Birkha\"{u}ser, Boston, 1985. 
    
\bibitem[GK]{GK} M. Gaberdiel and C. Keller, Differential operators for elliptic genera, arXiv: 0904.1831v1[hep-th], April 2009.
    
    \bibitem[H]{H} G. H\"{o}hn, Selbstdual Vertexoperatorsuperalgebren und das Babymonster, Bonner Math. Shriften \textbf{286}, 1996, vi+85pp.
  
 \bibitem[K]{K} V. Kac, \textit{Infinite-dimensional Lie algebras, 3rd. ed.}, C.U.P, 1990.

  
  \bibitem[KP]{KP} V. Kac and D. Peterson, Infinite-Dimensional Lie Algebras, Theta functions and Modular Forms, Adv. Math. \textbf{53} (1984), 125-264.
  
  \bibitem[KYY]{KYY} T. Kawai, Y. Yamada and S. Yang, Elliptic genera and $N=2$ superconformal field theory,
  Nucl. Phys. B \textbf{414} (1994), 191-
  
  
  \bibitem[L]{L} \textit{Elliptic Curves and Modular Forms in Algebraic Topology}
  (P. Landweber Ed.), Lect. Notes in Math. \textbf{1326}, Springer, Berlin, 1988.
  
  \bibitem[L1]{L1} H. Li, Symmetric invariant bilinear forms on vertex operator algebras,
  J. Pure  Appl. Math \textbf{96} (1994), 279-297.
  
  \bibitem[L2]{L2} H. Li, Local systems of twisted vertex operators, vertex operator superalgebras and twisted modules,
in \textit{Moonshine, the Monster, and related topics} Contemp. Math. Vol. \textbf{193} (1996), 203-236.
  
  \bibitem[LL]{LL} J. Lepowsky and H. Li, \textit{Introduction to Vertex Operator Algebras
  and Their Representations}, Birkha\"{u}ser, Boston, 2004.
 
 \bibitem[Ma]{Ma} G. Mason, Rational Vertex Operator Algebras and their Orbifolds, in 
 \textit{Moonshine: The First Quarter Century and Beyond}, Lond. Math. Soc. Lecture Note Series \textbf{372},
 J. Lepowsky, J. Mckay and M. Tuite Eds., C.U.P., 2010, 270-281.
 
  
  \bibitem[M]{M} M. Miyamoto, A modular invariance on the theta functions defined on vertex operator algebras, Duke Math. J. \textbf{101} No. 2 (2000), 221- 236.
  

\bibitem[EOT]{EOT} T. Eguchi, H. Ooguri and Y. Tachikawa, Notes on the K3 surface
and the Mathieu group M24, arXiv:1004.0956 [hep-th].
 
 \bibitem[Z]{Z} Y. Zhu, Modular-invariance of characters of vertex operator algebras,
 JAMS. \textbf{9}, No.1 (1996), 237--302.
 
 \end{thebibliography}
\end{document}